\documentclass{article}
\usepackage{amsmath,amssymb,amsthm} 
\input rlepsf

\oddsidemargin = \evensidemargin
\def\co{\colon\thinspace}
\newcommand{\Int}{\mbox{\rm Int}}
\newcommand{\SL}{\mbox{\rm SL}}

\newcommand{\op}{\overline{p}}
\newcommand{\oK}{\overline{K}}

\newtheorem{thm}{Theorem}
\newtheorem{lem}[thm]{Lemma}
\newtheorem{prop}[thm]{Proposition}
\newtheorem{cor}[thm]{Corollary}

\begin{document}

\title{A Legendrian surgery presentation of\\contact $3$--manifolds}
\author{Fan Ding and Hansj\"org Geiges}
\date{}

\maketitle

\begin{abstract}
\noindent
We prove
that every closed, connected contact $3$--manifold can be obtained
from $S^3$ with its standard contact structure by contact $(\pm 1)$--surgery
along a Legendrian link. As a corollary, we derive a result of Etnyre
and Honda about symplectic cobordisms (in slightly stronger form).
\end{abstract}

\vspace{2mm}

\noindent {\bf Mathematics Subject Classification (2000).}
Primary 53D35; Secondary 57R65, 57R90.

\vspace{2mm}

\noindent {\bf Keywords.} Contact structure, Legendrian link, contact
surgery, symplectic cobordism.


\section{Introduction}

A {\it contact structure} $\xi$ on a differential $3$--manifold $M$ is
a smooth tangent $2$--plane field such that any differential $1$--form
$\alpha$ that locally defines $\xi$ as $\xi =\ker\alpha$ satisfies the
nonintegrability condition $\alpha\wedge d\alpha\neq 0$. Notice that the
sign of $\alpha\wedge d\alpha$ only depends on~$\xi$, not the choice
of~$\alpha$. In particular, $M$ has to be orientable. Without any
essential loss of generality, we assume our $3$--manifolds to be
{\it oriented}, and we restrict attention to {\it positive} contact
structures, defined by the requirement that $\alpha\wedge d\alpha$
be a {\it positive} volume form. Moreover,
our contact structures will be understood to be {\it coorientable},
which is equivalent to saying that a defining $1$--form $\alpha$ exists
globally.

The {\it standard contact structure} $\xi_0$ on $S^3\subset {\mathbb R}^4$
is defined by $\sum_{i=1}^2(x_idy_i-y_idx_i)=0$. A theorem of
Lutz~\cite{lutz70} and Martinet~\cite{mart71} states that on any given
(closed, orientable) $3$--manifold $M$ one can find a contact structure
in each homotopy class of $2$--plane fields by performing surgery
on $(S^3,\xi_0)$ along a link {\it transverse} to~$\xi_0$. (The part
about `each homotopy class of $2$--plane fields' is not completely covered
by the cited references, but belongs to contact geometric folklore,
cf.~\cite{geig}.)

Recall the dichotomy between {\it overtwisted} and {\it tight} contact
structures. A contact structure $\xi$ on $M$ is called overtwisted
if one can find an embedded disc $D\hookrightarrow M$ with boundary
$\partial D$ tangent to~$\xi$, but $D$ transverse to $\xi$ along~$\partial D$.
If no such $D$ exists, then $\xi$ is called tight. Eliashberg has
shown that the (isotopy) classification of
overtwisted contact structures coincides with the
(homotopy) classification
of tangent $2$--plane fields~\cite{elia89}.

In the Lutz-Martinet theorem no statement is made whether {\it every}
contact structure on $M$ can be obtained via transverse contact surgery
as described, but from the work of Eliashberg~\cite{elia89}
one can deduce that
at least all overtwisted contact structures (up to isotopy)
are covered by this construction.

A smooth knot $K$ in a contact manifold $(M,\xi )$ is called {\it Legendrian}
if it is everywhere tangent to~$\xi$. Such a Legendrian knot inherits a
canonical {\it contact framing} of its normal bundle, defined via
any vector field in $\xi |K$ transverse to~$K$. In~\cite{dige01}
we described a notion of contact $r$--surgery along such Legendrian
knots, where the surgery coefficient $r\in {\mathbf Q}\cup\{
\infty\}$ is measured with respect to the contact framing; details will be
recalled below. In general the resulting contact structure on the
surgered manifold depends on choices, but for $r=1/k$ with $k\in
{\mathbb Z}$ it is fully determined by $(M,\xi )$, $K$
and~$r$~\cite[Prop.~7]{dige01} (up to contactomorphism,
i.e.\ diffeomorphism preserving the contact structures).
The inverse of a contact $(1/k)$--surgery is a contact $(-1/k)$--surgery;
for the case $k=\pm 1$ (which is the one most relevant to our discussion)
we give an explicit proof of this fact in Section~\ref{section:pm1}.

Our main theorem is the following.

\begin{thm}
\label{thm:main}
Let $(M,\xi )$ be a closed, connected
contact $3$--manifold. Then $(M,\xi )$
can be obtained by contact $(\pm 1)$--surgery along a Legendrian
link in $(S^3,\xi_0)$.
\end{thm}

Contact $(-1)$--surgery coincides with the symplectic handlebody surgery
described by Eliashberg~\cite{elia90} and Weinstein~\cite{wein91},
cf.~\cite{gomp98} and Section~\ref{section:pm1} below.
In particular, if $(M',\xi ')$ is obtained from $(M,\xi )$ by contact
$(-1)$--surgery along a Legendrian link, then there is a so-called
{\it Stein cobordism} from $(M,\xi )$ to $(M',\xi ')$, see~\cite{etho}
for details, and~\cite{etny98} for background information.
This is denoted by $(M,\xi )\prec (M',\xi ')$. If
$(S^3,\xi_0)\prec (M',\xi ')$, then $(M',\xi ')$ is {\it Stein fillable};
the concave end $(S^3,\xi_0)$ of the Stein cobordism can be filled
with the standard symplectic $4$--ball.

It would be possible to deduce Theorem~\ref{thm:main} directly from
the results of Eliashberg~\cite{elia90}. By two contact $(+1)$--surgeries
on a given $(M,\xi )$ one can effect a simple Lutz twist
(Proposition~\ref{prop:lutz} of the present paper). This produces
an overtwisted contact structure on~$M$. Surgery along a suitable link in $M$
produces~$S^3$. Eliashberg's results, in particular~\cite[Lemma~2.4.3]{elia90},
then say that $S^3$ (with some contact structure determined by
the choice of surgeries)
can be obtained from $M$ by contact $(-1)$-surgeries. (Alternatively,
one may appeal
to \cite[Lemma~4.4]{gomp98} and then directly to~\cite[Thm.~1.3.4]{elia90}.)
As shown in Section~\ref{section:lutz} below, a further sequence
of contact $(-1)$--surgeries on this $S^3$ will produce $(S^3,\xi_0)$.

A.~Stipsicz~\cite{stip} has suggested a similar shortcut to the
proof of our theorem.

The extra benefit of the proof of Theorem~\ref{thm:main} given in the present
paper lies in the fact that it contains an algorithm for turning
contact surgeries with rational coefficients into surgeries of
the desired kind.  It thus allows to discuss explicit surgery descriptions
of contact manifolds. We hope to elaborate on this point in a future
publication.

It is now an easy matter to derive from Theorem~\ref{thm:main}
the following corollary, which
generalises a theorem of Etnyre and Honda~\cite[Thm.~1.1]{etho}.

\begin{cor}
Let $(M_i,\xi_i)$ be contact $3$--manifolds (as in
Theorem~\ref{thm:main}), $i=1,\ldots n$. Then there
is a Stein fillable contact manifold $(M,\xi )$ such that $(M_i,\xi_i)\prec
(M,\xi )$ for each~$i$.
\end{cor}

\noindent {\bf Proof.}
Let $(M_i,\xi_i)$ be obtained from $(S^3,\xi_0)$ by contact
$(+1)$--surgery along a link $L_i^+$ and $(-1)$--surgery
along a link~$L_i^-$. We may assume that the links $L_1^+,L_1^-,\ldots ,
L_n^+,L_n^-$ are pairwise disjoint. Let $(M,\xi )$ be the contact manifold
obtained from $(S^3,\xi_0)$ by contact $(-1)$--surgery along the link
$L_1^-\cup\ldots\cup L_n^-$. Then $(M,\xi )$ is Stein fillable. Moreover,
each $(M_i,\xi_i)$ is obtained from $(M,\xi)$ by contact $(+1)$--surgeries,
so $(M,\xi )$ is obtained from $(M_i,\xi_i)$ by contact $(-1)$--surgeries.

\vspace{2mm}

Etnyre and Honda, by contrast, base their proof on a result of
Giroux~\cite{giro02},~\cite{giro}
about open book decompositions adapted to
contact structures. In principle, their methods suffice to deduce this
stronger corollary, too. Our proof arguably has the advantage that
a surgery description is more appropriate in the context of
a cobordism theoretic result than an open book decomposition.
\section{Contact surgery}
\label{section:surgery}
We recall a few definitions and results from~\cite{dige01}, chiefly to
fix notation and conventions.

Let $K$ be a Legendrian knot in a contact $3$--manifold~$(M,\xi )$.
{\it Rational surgery} on $K$ with coefficient $r=p/q\in {\mathbf Q}\cup
\{ \infty\}$ (with $p,q$ coprime) is defined as follows: Denote a tubular
neighbourhood of K
(diffeomorphic to a solid torus) by $\nu K$. Let $(\mu,\lambda)$ be a
positively oriented basis for $H_1(\partial\nu K;{\mathbf Z})\cong
{\mathbf Z}\oplus {\mathbf Z}$, where $\lambda$ is determined up to sign as
the class of a parallel copy of $K$ determined by the contact framing, and
$\mu$ is determined  by a suitably oriented meridian (i.e.\ a
nullhomologous circle in $\nu K$), cf.~\cite[p.~672]{gomp98}.
We obtain a new manifold $M'$ by cutting
$\nu K$ out of $M$ and regluing it by a diffeomorphism of $\partial(\nu K)$
sending $\mu$ (on the boundary of the solid torus to be glued in)
to  the curve $p\mu+q\lambda$ (on the boundary of $M-\nu K$).
This procedure determines $M'$ up to orientation-preserving diffeomorphism.
For $r=\infty$ the surgery is trivial.

There is a unique contact geometric model for the tubular neighbourhood 
of a Le\-gendrian knot. In order to describe this, we
consider $N={\mathbf R}^2\times ({\mathbf R}/{\mathbf Z})$ with coordinates
$(x,y,z)$ and contact structure $\zeta$ defined by
\[ \cos (2\pi z)\, dx-\sin (2\pi z)\, dy =0.\]
For each $\delta >0$, let
\[ N_{\delta}= \{ (x,y,z)\in N\co x^2+y^2\le {\delta}^2\}.\]
We identify $T_{\delta}=\partial N_{\delta}$ with ${\mathbf R}^2/{\mathbf Z}^2$,
using the contact framing, and
write $(\mu ,\lambda )$ for a positively oriented basis for
$H_1(T_{\delta};{\mathbf Z})\cong {\mathbf Z}\oplus {\mathbf Z}$,
with $\mu$ corresponding to a meridian and $\lambda$ to a longitude
determined by this framing. A possible representative of $\lambda$ would be
\[ \{ (\delta\sin (2\pi z),\delta\cos (2\pi z),z)\co z\in
{\mathbf R}/{\mathbf Z}\}.\]

The vector field $X=x\frac{\partial}{\partial x}+y
\frac{\partial}{\partial y}$ is a {\it contact vector field} for the
contact structure~$\zeta$, i.e.\ its flow preserves $\zeta$. This
vector field is
transverse to~$T_{\delta}$. Any closed surface in a contact manifold
that admits a transverse contact vector field is called
{\it convex}~\cite{giro91}.

The set $\Gamma_{T_{\delta}}$ of points $w\in T_{\delta}$
where $X(w)\in\zeta (w)$ is equal to
\[ \Gamma_{T_{\delta}}=\{ (\pm\delta\sin (2\pi z),\pm\delta\cos
(2\pi z),z)\co z\in {\mathbf R}/{\mathbf Z}\} .\]
This set is called the {\it dividing set} of $T_{\delta}$; its
isotopy type is independent of the choice of contact vector field
transverse to $T_{\delta}$. Notice that the dividing set is made up
of two copies of the longitude determined by the contact framing.

In general, if $T$ is a convex torus in a tight contact $3$--manifold,
the dividing set $\Gamma_T$ consists of an even number $\#\Gamma_T$
of parallel essential curves. (Conversely, the absence of homotopically
trivial curves on a convex surface $\Sigma\neq S^2$ in a contact
manifold guarantees that $\Sigma$ has a tight neighbourhood.) This is
(a special case of) Giroux's criterion, cf.~\cite[Thm.~3.5]{hond00}. After
a diffeomorphism isotopic to the identity, one may assume the dividing
curves to be linear relative to an identification of $T$ with
${\mathbf R}^2/{\mathbf Z}^2$. Their slope will be called the {\it slope}
$s(T)$ of the convex torus~$T$. If $T$ bounds a solid torus, our convention will
be to identify it with ${\mathbf R}^2/{\mathbf Z}^2$ in such a way that
the meridian $\mu$ has slope~$0$ and the (chosen) longitude $\lambda$
has slope~$\infty$. Thus, with respect to the contact framing,
the $T_{\delta}$ above is a convex torus with $\#\Gamma_{T_{\delta}}=2$
and~$s(T_{\delta})=\infty$. The key observation by Giroux is that
the dividing set of a convex surface $\Sigma$ encodes all the essential
contact geometric information in a neighbourhood of~$\Sigma$. In particular,
one can glue contact manifolds along surfaces with the same dividing set
(taking boundary orientations into account). The specific results of Giroux
necessary to make this statement precise are summarised in Section~2
of~\cite{dige01}.

We now describe how to perform this rational surgery on $K$ in such a
way that the resulting manifold again carries a contact structure. Write
\[ C=\{ (x,y,z)\in N\co x=y=0\} \]
for the spine of~$N$. Then there is a contact embedding $f\co (N_2,\zeta)\to
(M,\xi)$ such that $f(C)=K$. We want to construct a contact structure
$\xi '$ on the manifold $M'$ obtained from $M$ by rational surgery on
$K$ with coefficient~$r=p/q\in {\mathbf Q}\cup \{\infty\}$. Let
\[ P=\{(x,y,z)\in N\co 1\le x^2+y^2\le 4\}=N_2\setminus\Int\, N_1.\]
Let $g\co P\to P$ be an orientation-preserving diffeomorphism
sending $T_{\delta}$ to $T_{\delta}$, ${\delta}=1,2$, and
$\mu$ to $p\mu +q\lambda$.

If $p\neq 0$, then
$(g_*)^{-1}(\zeta)$ is a contact structure on $P$ with respect to
which $T_{\delta}$ is a convex torus of
non-zero slope. By~\cite[Thm.~2.3]{hond00}, which gives
an enumeration of tight contact structures on the solid torus with
convex boundary as in our situation (and in particular shows this
set of contact structures to be non-empty), the contact structure
$(g_*)^{-1}(\zeta )$ on $P$ can be extended to a tight
contact structure $\zeta '$ on~$N_2$. Define
\[ M' =(M-f(N_1))\cup N_2/\sim ,\]
where $w\in P\subset N_2$ is identified with $f(g(w))\in M$,
and equip $M'$ with the contact structure $\xi '$ it inherits
from $(M,\xi )$ and $(N_2,\zeta ')$. We say that $(M',\xi ')$ is obtained from
$(M,\xi )$ by {\it contact $r$--surgery on}~$K$, $r\neq 0$.

Topologically, $M'$ is completely determined by $K$ and~$r$, but $\xi '$
depends on the choice of~$\zeta '$. Only for $r=1/k$, $k\in {\mathbf Z}$,
this choice is unique. Moreover, as pointed out in the introduction,
contact $(-1/k)$--surgery is the inverse of contact
$(1/k)$--surgery~\cite[Prop.~8]{dige01}; see the following section for
the case~$k=\pm 1$.

If $p=0$ (and $q=\pm 1$), then $(g_*)^{-1}(\zeta)$ will be a contact
structure on $P$ with respect to
which $T_{\delta}$ is a convex torus of slope zero. In that case we use
a construction that will appear again in various guises in the discussion
of Lutz twists below. Namely, consider a solid torus $S^1\times D^2$
with $S^1$--coordinate $\theta$ and polar coordinates $r,\varphi$
on~$D^2$. Write $\mu_0$ for the meridian of this solid torus and
$\lambda_0$ for a longitude given by setting $\varphi$ equal to a constant.
Let $\zeta_0$ be the contact structure on $S^1\times D^2$ defined by
the differential $1$--form
$\beta_0 =h_1(r)\, d\theta +h_2(r)\, d\varphi$. Here we impose the following
general conditions on the smooth plane curve $r\mapsto
\gamma (r)=(h_1(r),h_2(r))$:
\begin{itemize}
\item $h_1(r)\equiv\pm 1$ and $h_2(r)=\pm r^2$ for $r$ near~$0$; this
guarantees that the $1$--form $\beta_0$
is actually defined at $r=0$.
\item Position vector $\gamma (r)$ and velocity vector $\gamma '(r)$ are
never parallel to each other (in particular $\gamma (r)\neq (0,0)$);
this ensures that $\beta_0=0$ does indeed
define a contact structure, cf.~\cite{geig01}.
\end{itemize}
Without loss of generality we shall always assume that $\gamma (r)$ winds
in counter-clockwise (i.e.\ positive) direction around the origin.
This amounts to
fixing the orientation of $S^1\times D^2$ as the one given by
$d\theta\wedge r\, dr\wedge d\varphi$.

In order to define $0$--surgery, we further require that $h_1(1)=1$ and
$h_2(1)=0$, and that $\gamma (r)$ completes exactly one half-turn
as $r$ increases from $0$ to~$1$. Then along the boundary $T^2=\partial
(S^1\times D^2)$ the contact structure is given by $d\theta =0$.
Recall that the {\it characteristic foliation} $\Sigma_{\xi}$
of a surface $\Sigma$
in a contact manifold $(M,\xi )$ is the singular $1$--dimensional
foliation defined by the intersection of the tangent planes of
$\Sigma$ with~$\xi$. Hence, the characteristic foliation $T^2_{\zeta_0}$
consists of the meridional loops. By a $C^{\infty}$--small perturbation
of $T^2$ (inside a slight thickening of $S^1\times D^2$) one can turn
this $2$--torus into a convex surface with $\#\Gamma_{T^2}=2$ and
slope~$0$, see~\cite[p.~795/6]{giro99}. Then the solid torus with this
convex boundary can be glued into $M-f(N_1)$, using the attaching map
corresponding to $p=0$, i.e.\ the map that sends $\mu_0$ to $\lambda$
and $\lambda_0$ to~$-\mu$, say. (Notice, in particular, that the proviso
`if such exists' in Remark~(1), Section~3 of~\cite{dige01} was superfluous.)
\section{Contact $(\pm 1)$--surgeries}
\label{section:pm1}
We follow Weinstein~\cite{wein91} in the description of contact
$(\pm 1)$--surgeries as symplectic handlebody surgeries. Consider
${\mathbb R}^4$ with coordinates $x,y,z,t$ and standard symplectic
form $\omega =dx\wedge dy+dz\wedge dt$. The vector field
\[ Z=2x\partial_x-y\partial_y+2z\partial_z-t\partial_t\]
is a Liouville vector field for~$\omega$, that is,
\[ L_Z\omega=d(i_Z\omega)=\omega.\]
This implies that $\alpha =i_Z\omega$ defines a contact form on
any hypersurface in ${\mathbb R}^4$ transverse to~$Z$.

First we consider the hypersurface $V=\{ y^2+t^2=\varepsilon^2 \}$.
Introduce a new coordinate $\theta$ by setting $y=\varepsilon\cos\theta$,
$t=\varepsilon\sin\theta$. Then
\[ \alpha |_V=\varepsilon\bigl( (-2x\,\sin\theta +2z\,\cos\theta )\, d\theta
+\cos\theta\, dx+\sin\theta\, dz\bigr) .\]
Observe that the circle $L= V\cap\{ x=z=0\}$ is Legendrian, and that
the radial vector field $x\partial_x+z\partial_z$ is a contact vector field.
The contact structure near $L$ is not quite in the normal form considered
in the previous section, but it is sufficient to know that small enough
neighbourhoods of any two Legendrian knots are contactomorphic.

\begin{figure}[h]
\centerline{\relabelbox\small
\epsfxsize 10cm \epsfbox{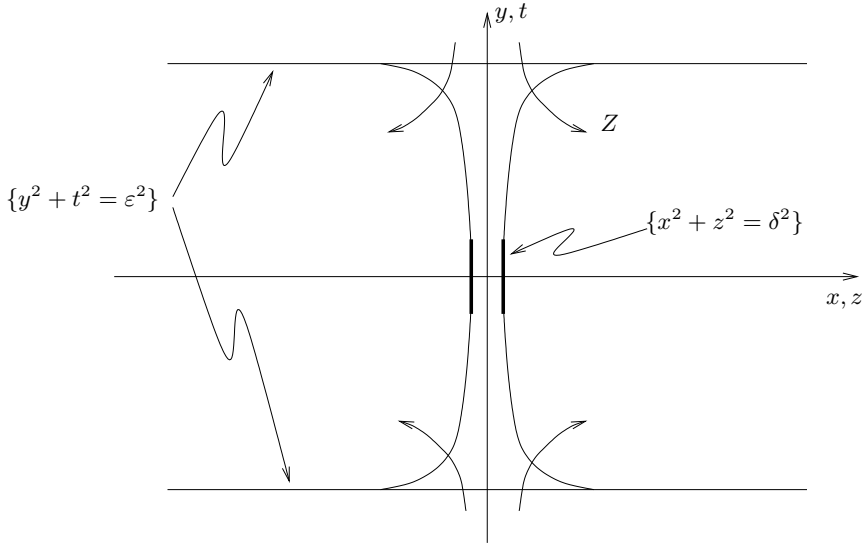}
\extralabel <-3.5cm, 5.5cm> {$Z$}
\extralabel <-.5cm, 3.2cm> {$x,z$}
\extralabel <-4.9cm, 7cm> {$y,t$}
\extralabel <-11.4cm, 4.5cm> {$\{ y^2+t^2=\varepsilon^2\}$}
\extralabel <-2.9cm, 4.2cm> {$\{ x^2+z^2=\delta^2\}$}
\endrelabelbox}
\caption{Contact $(-1)$--surgery.\label{figure:min1surgery}}
\end{figure}

Figure~\ref{figure:min1surgery} shows how to attach a $2$--handle (with
boundary transverse to~$Z$) along a tubular neighbourhood of $L$ in~$V$.
We want to show that this amounts to a contact $(-1)$--surgery along~$L$:
The orientation on $V$ is given by $\alpha\wedge d\alpha =
\varepsilon^2 dx\wedge dz\wedge
d\theta$, and the contact framing of $L$ in $V$ by $\sin\theta\,\partial_x
-\cos\theta\,\partial_z$. This implies that the framing of the surgery,
given by~$\partial_x$, makes one negative twist in the $xz$--plane
relative to the contact framing as we go once around $L$ in
$\theta$--direction (this sign is independent of the choice of
orientation for~$L$).

Given a Legendrian knot $K$ in a contact manifold $(M,\xi )$,
a small neighbourhood of
$K$ is contactomorphic to a neighbourhood of $L$ in~$V$. Therefore the
handle described in Figure~\ref{figure:min1surgery} can be used to perform
a contact $(-1)$--surgery on~$K$.

If $(M,\xi )$ is strongly symplectically fillable ($\emptyset\prec (M,\xi )$
in the notation of the introduction), that is, $M$ is the boundary of
a symplectic manifold $(W,\omega )$ admitting a Liouville vector field
$Z$ near $\partial W=M$, transverse to $M$, pointing outwards, and
satisfying $\ker i_Z\omega =\xi$, then the construction shows
that the surgered manifold will again be strongly symplectically
fillable (for this one also needs to `glue' the Liouville vector field
on $W$ near $K\subset M$ and that in the handle picture near $L$,
see~\cite{wein91} for the technical details).

\begin{figure}[h]
\centerline{\relabelbox\small
\epsfxsize 10cm \epsfbox{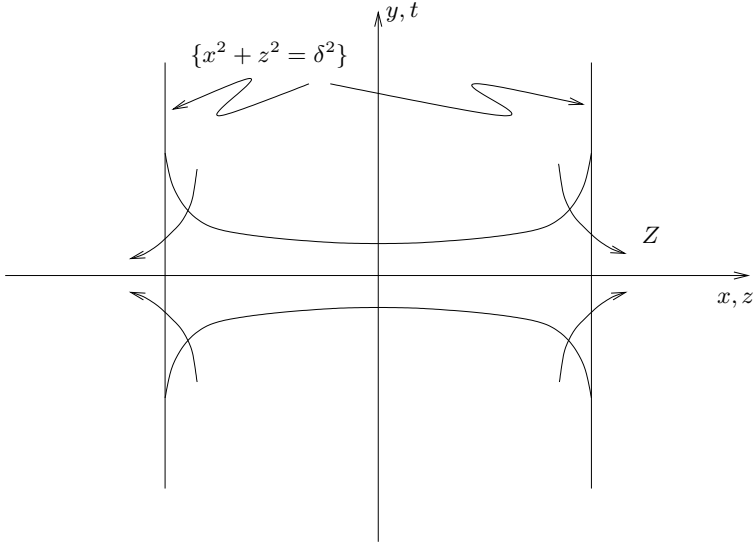}
\extralabel <-.5cm, 3.2cm> {$x,z$}
\extralabel <-4.9cm, 7cm> {$y,t$}
\extralabel <-1.5cm, 4cm> {$Z$}
\extralabel <-7.5cm, 6.4cm> {$\{ x^2+z^2=\delta^2\}$}
\endrelabelbox}
\caption{Contact $(+1)$--surgery.\label{figure:1surgery}}
\end{figure}

On $V'=\{ x^2+z^2=\delta^2\}$ by contrast, and with $x=\delta\cos\varphi$,
$z=\delta\sin\varphi$, we have
\[ \alpha |_{V'}=\delta\bigl( (-y\sin\varphi\, d\varphi +t\cos\varphi\,
d\varphi )+ 2\cos\varphi\, dy +2\sin\varphi\, dt\bigr) .\]
Now we have $\alpha\wedge d\alpha = 2\delta^2 dt\wedge dy\wedge d\varphi$,
and the contact framing of $L'=V'\cap\{ y=t=0\}$ is given by
$\sin\varphi\,\partial_y-\cos\varphi\,\partial_t$. So the framing $\partial_y$
of the surgery indicated in Figure~\ref{figure:1surgery} makes one positive
twist in the $ty$--plane relative to the contact framing when going once
around $L'$ in $\varphi$--direction. Hence this amounts to a contact
$(+1)$--surgery. It can be performed as a symplectic handlebody surgery
on the concave end of a symplectic cobordism (where the Liouville vector field
points inward).

\begin{figure}[h]
\centerline{\relabelbox\small
\epsfxsize 6cm \epsfbox{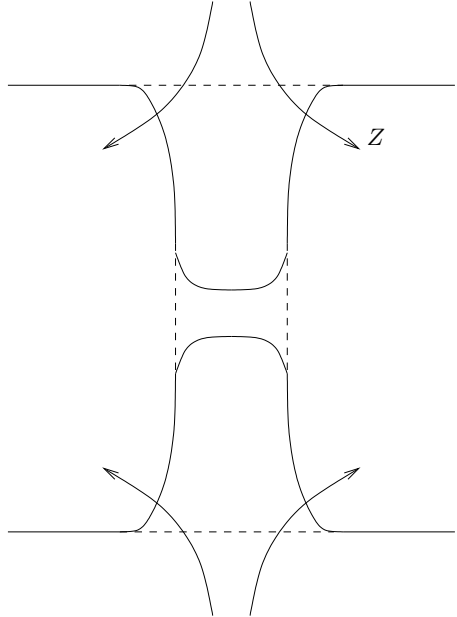}
\extralabel <-1.2cm, 6.3cm> {$Z$}
\endrelabelbox}
\caption{Composition of a contact $(-1)$--surgery and a $(+1)$--surgery.
\label{figure:compose}}
\end{figure}

Figure~\ref{figure:compose} shows that the composition of a contact
$(-1)$--surgery and a contact $(+1)$--surgery leads to a manifold
contactomorphic to the one we started from. 
Given two hypersurfaces in a symplectic manifold, both transverse to the same
Liouville vector field $Z$ and with contact structures induced by~$Z$,
one can show by a straightforward calculation
(cf.~\cite[Lemma~2.2]{wein91}) that the local diffeomorphism
between the two hypersurfaces defined by following the integral curves of $Z$
is a contactomorphism. Thus, we simply have to move the `bump' in
Figure~\ref{figure:compose} along the integral curves
of $Z$ back to the original hypersurface~$V$.

The same figure gives an elementary proof (not appealing to
an $h$--principle, for instance)
that an iteration of contact $(\pm 1)$--surgeries may
be replaced by a simultaneous contact surgery (without changing the
surgery coefficients) along a suitable link. This amounts,
as it were, to a handle reordering lemma for
contact surgeries. Suppose we have performed a surgery along a Legendrian
knot $K_1$ as in Figure~\ref{figure:min1surgery}, and we now want to
carry out a further surgery along a Legendrian knot $K_2$ in the surgered
manifold. By a contact isotopy we can move $K_2$ away from the belt
sphere of the first surgery (such contact isotopies can easily be constructed
using contact Hamiltonians, cf.~\cite{geig}). As before, we see that the
complement of a neighbourhood of this belt sphere is
contactomorphic to the complement of a neighbourhood of~$K_1$ in the
original manifold. Using the radial contact vector field near~$K_1$,
we can move $K_2$ away from the first handle.
\section{Front projections}
As a first preparation for the proof of Theorem~\ref{thm:main} we recall
how to describe Legendrian links in ${\mathbf R}^3$ with its standard tight
contact structure $dz+x\, dy=0$ (which is contactomorphic to $(S^3,\xi_0)$
with a point removed, and the universal local model for contact structures).

We visualise links in ${\mathbf R}^3$ by projecting them into the
$(y,z)$--plane.
The condition for $p(t)=(x(t),y(t),z(t))$ to describe a Legendrian curve in
${\mathbf R}^3$ is that $z'(t)+x(t)y'(t)=0$. Any closed curve $\op (t)
=(y(t),z(t))$ in the $(y,z)$--plane that may have cusps and transverse
self-crossings, but does not have any vertical tangencies (such a curve
will be called a {\it front}), comes from
a unique Legendrian knot $p(t)$ in ${\mathbf R}^3$, found by setting
$-x(t)$ equal to the slope of~$\op (t)$. Cusps of $\op (t)$ correspond
to points where $p(t)$ is parallel to the $x$--axis.

The plane curve $\op (t)$ is called the {\it front projection}
of~$p(t)$. The theory of front projections was developed by
Arnol'd~\cite{arno80}; all information relevant for our purposes can be found
in Section~1 of~\cite{gomp98}.

There are two Legendrian isotopy invariants of oriented Legendrian knots~$K$,
viz.\ the Thurston-Bennequin invariant $tb(K)$ and the rotation number~$r(K)$.
The Thurston-Bennequin invariant measures the linking number of $K$ with
its push-off determined by the contact framing;
we refer to~\cite{gomp98} for the definition of the rotation number.
In the present context we
only need to know the following properties of these invariants:

\begin{itemize}
\item If the front projection $\oK$ of the Legendrian knot $K\subset
{\mathbf R}^3$ has no crossings, then $tb(K)=-c(\oK )/2$, where $c(\oK)$
denotes the number of cusps of the front projection. (This is immediate from
the definition of~$tb(K)$.)
\item The rotation number $r(K)$ is equal to $\lambda_-(\oK )-\rho_+
(\oK )$, where $\lambda_-$ denotes the number of cusps with vertex on
the left and oriented downwards, $\rho_+$ the number of right cusps
oriented upwards (remember that $K$ has to be oriented). Changing the
orientation of $K$ will change the sign of~$r(K)$ (and leave $tb(K)$
unaltered).
\item For the trivial knot~$K$,
these invariants satisfy the Thurston-Bennequin inequality
\[ tb(K)+|r(K)|\leq -1.\]
\end{itemize}
The consequences of these properties that we shall use in the next section
are the following:
\begin{itemize}
\item[(i)] The Legendrian knot $K_0$ in ${\mathbf R}^3$ whose front projection
is a simple closed curve with two cusps has the invariants $tb(K_0)=-1$
and~$r(K_0)=0$.
\item[(ii)] By approximating the front projection of $K_0$ by a front with
an additional zigzag (a left and a right cusp, oriented up or down)
and no self-crossings, we can decrease the Thurston-Bennequin
invariant by~$1$ (i.e.\ add a negative twist to the contact framing)
and decrease or increase the rotation number by~$1$. Hence, by adding $n$
zigzags we can realise $tb=-n-1$ and any of the $n+1$ different rotation
numbers
\[ -n,-n+2,-n+4,\ldots ,n-2,n,\]
i.e.\ all values in steps by~$2$ in the range allowed by the
Thurston-Bennequin inequality.
\end{itemize}
Another useful fact that can easily be proved using the concept of
front projections is that any knot $p(t)=(x(t),y(t),z(t))$ in
${\mathbf R}^3$ can be $C^0$--approximated by a Legendrian knot: simply
approximate the plane curve $(y(t),z(t))$ by a front (with many zigzags)
whose slope approximates~$-x(t)$. The same is true for knots in arbitrary
contact manifolds.
\section{Equivalent contact surgeries}
The strategy for proving Theorem~\ref{thm:main} is as follows. Let
$(M,\xi )$ be given. By the topological surgery presentation theorem
of Lickorish~\cite{lick62} and Wallace~\cite{wall60}, the manifold $M$
can be obtained by performing integer surgery along a link in~$S^3$.
Thus, conversely, we can recover $S^3$ by integer surgery along a link in~$M$.
By the remark at the end of the preceding section, we may assume this to
be a Legendrian link, so that we can perform the surgeries as contact integer
surgeries. The final product of these surgeries will be $S^3$ with some
contact structure~$\xi_0'$.

In the present section we show that all these contact integer surgeries
may be replaced by contact $(\pm 1)$--surgeries along a different link.
In the next section we show how to transform $(S^3,\xi_0')$ into
$(S^3,\xi_0)$ by further contact $(\pm 1)$--surgeries. Notice that the
surgery curves may be moved about by contact isotopies (which exist
in abundance due to their constructibility via contact Hamiltonians);
this allows to replace an iteration of contact surgeries by a simultaneous
contact surgery (without altering the surgery coefficients) along a
suitable link.

\begin{prop}
\label{prop:r<0}
Any contact $r$--surgery with $r<0$ can be obtained by a sequence of contact
$(-1)$--surgeries.
\end{prop}

\noindent {\bf Proof.}
We want to show that a contact $r$--surgery, $r<0$, along a Legendrian
knot $K$ may be replaced by a sequence of contact $(-1)$--surgeries
along suitable Legendrian knots inside a tubular neighbourhood
$\nu K$ of~$K$. This would seem to suggest that one
should identify $\nu K$ with $S^1\times D^2$ using the contact framing
of~$K$. For our inductive procedure, however, it turns out to be more
opportune to twist this framing by~$-1$.

Thus, we set $N_0=S^1\times D^2$ and write $\mu_0$ for its meridian
and $\lambda_0$ for the longitude determined by setting $\varphi$
equal to a constant (in the usual coordinates $(\theta ,r,\varphi )$
on $S^1\times D^2$). Recall the following result~\cite[Thm.~8.2]{kand97},
cf.~\cite[Prop.~4.3]{hond00} (which is the reason why contact
$(1/k)$--surgery is uniquely defined, see~\cite[Prop.~7]{dige01}):

\begin{prop}
\label{prop:1/k}
For any integer $k$ (including~$0$) there exists a unique (up to
isotopy fixed at the boundary) tight contact structure on $S^1\times D^2$
with a fixed convex boundary with $\#\Gamma_{\partial (S^1\times D^2)}=2$
and slope $s(\partial (S^1\times D^2))=1/k$.
\end{prop}

Let $\xi$ be the essentially unique tight contact structure on $N_0$
with $\#\Gamma_{\partial N_0}=2$ and boundary slope~$-1$. In other words,
the longitude determined by the contact framing is $\lambda_0-\mu_0$.
Thus, if $\mu$ denotes the meridian of $\nu K$
and $\lambda$ the longitude of $\nu K$ determined by the contact framing,
then the identification of $N_0$ with $\nu K$ is given by
\[ \mu_0\longmapsto \mu,\;\; \lambda_0\longmapsto \mu+\lambda,\]
so that $\lambda_0-\mu_0\mapsto \lambda$.

If we let the meridian correspond to the first coordinate direction
in ${\mathbf R}^2/{\mathbf Z}^2$, and the longitude to the second, then
the above identification is described by the matrix
$\left(\begin{array}{cc}1&1\\0&1\end{array}\right)$.

Now let $K_1$ be a Legendrian knot in $N_0$, isotopic to the spine
$S^1\times \{ 0\}$, and with boundary slope $1/(r_1+1)$,
$r_1+1\leq -1$. That is, $K_1$ has a tubular neighbourhood $N_1\subset N_0$
(with convex boundary $\partial N_1$ with $\#\Gamma_{\partial N_1}=2$)
such that, under the diffeomorphism of $\partial N_1$ with $\partial N_0$
determined by the isotopy of $K_1$ with $S^1\times \{ 0\}$, the
longitude of $\partial N_1$ determined by the contact framing
is equal to $\lambda_0+(r_1+1)\mu_0$ (and meridian still equal to~$\mu_0$).
For $r_1+1\leq -1$ this is exactly what remark (ii) in the preceding
section allows us to do.

The subscript $i=0,1$ in $N_i$ here is of course just a counter and does not
refer to the radius of these solid tori. The same applies to all subsequent
occurences of the notation `$N_i$'. There should be little grounds for
confusion with the notation in Section~\ref{section:surgery} in the
definition of contact $r$--surgery.

Now perform a contact $(-1)$--surgery on~$K_1$. That is, we cut out $N_1$
and glue back a solid torus $N_1'$ by sending its meridian $\mu_1$ to
\[ \mu_0-(\lambda_0+(r_1+1)\mu_0)=-r_1\mu_0-\lambda_0.\]
We may choose a longitude $\lambda_1$ for $N_1'$ such that the
gluing is described by the matrix $\left(\begin{array}{cc}-r_1&1\\-1&0
\end{array}\right)\in\SL (2,{\mathbf Z})$
with respect to $(\mu_1,\lambda_1)$ and
$(\mu_0,\lambda_0)$ (that is, $\lambda_1$ maps to $\mu_0$). Notice that
the surgered $N_0$ is still a solid torus: after cutting out $N_1$
we have a torus shell $T^2\times [0,1]$, and we glue a solid torus
$N_1'$ to one of its boundary components.

Observe that the curve $\lambda_1-\mu_1$ on $\partial N_1'$ is mapped to
a dividing curve
\[ \lambda_0+(r_1+1)\mu_0\]
on $\partial N_1$. This
determines the extension of the contact structure over~$N_1'$;
it is the unique tight contact structure for which $\partial N_1'$
is convex with $\#\Gamma_{\partial N_1'}=2$ and $s(\partial N_1')=-1$.

We now let $N_1'$ take the role of~$N_0$. That is, we choose a Legendrian
knot isotopic to the spine of~$N_1'$ and with boundary slope
$1/(r_2+1)$ (expressed in terms of $(\mu_1,\lambda_1)$), where
$r_2+1\leq -1$, and we perform contact $(-1)$--surgery on this knot.
Continuing this process, we see that by a sequence of
$(-1)$--surgeries (say $n$ of them) we can realise a topological $r$--surgery,
where $r=p/q$ is related to the integers $r_1,\ldots ,r_n\leq -2$ by
\[ \left(\begin{array}{cc}1&1\\0&1\end{array}\right)
   \left(\begin{array}{cc}-r_1&1\\-1&0\end{array}\right)
   \left(\begin{array}{cc}-r_2&1\\-1&0\end{array}\right)
   \cdots
   \left(\begin{array}{cc}-r_n&1\\-1&0\end{array}\right)
   =
   \left(\begin{array}{cc}p&p'\\q&q'\end{array}\right) .\]
From this one finds that $p/q$ has the continued fraction expansion
\[ r_1+1-\cfrac{1}{r_2-\cfrac{1}{\dotsb -\cfrac{1}{r_n} }}\; ,\]
see~\cite[Thm.~7.1.1]{rose88}, for instance. We abbreviate this
continued fraction as $[r_1+1,r_2,\ldots ,r_n]$. The following lemma we
leave as a simple exercise (use the Euclidean algorithm with
negative remainders).

\begin{lem}
Any rational number $r<0$ can be written in the form
\[ r=[r_1+1,r_2, \ldots ,r_n] \]
with $r_1,\ldots ,r_n\leq -2$.
\end{lem}

The next thing to notice is that the topological $r$--surgery we obtain in
this way is actually a contact $r$--surgery. For this we only need to
show that we still have a tight contact structure on the surgered~$N_0$
(in general, although contact $r$--surgery is defined by choosing a tight
extension over the solid torus to be glued in, the resulting contact
structure on the surgered manifold need not be tight). But this
follows from the fact that contact $(-1)$--surgeries correspond to
symplectic handlebody surgeries. Furthermore, $N_0$ has a contact
embedding into $(S^3,\xi_0)$, which is symplectically filled by
the standard $4$--ball~$B^4$. So we can realise the surgered $N_0$ as a subset
of the manifold obtained as the new boundary after attaching $n$
symplectic $2$--handles to $B^4$ along $S^3=\partial B^4$. As a symplectically
fillable contact manifold, it is tight. See~\cite{etny98} for more
details on the notions used in this paragraph.

It remains to show that {\it any} contact $r$--surgery can be performed in
this way. The net result of performing the sequence of $(-1)$--surgeries
inside $N_0$ is a solid torus with a tight contact structure and with
convex boundary $\partial N_0$ whose dividing set has two
components (for the boundary has remained unchanged). The slope
of the dividing set has changed, however. To determine this slope,
we need to find the curve on $\partial N_n'$ (the boundary of the 
last solid torus to be glued in) that maps to $\lambda_0-\mu_0$
under the successive gluings. That is, $s(\partial N_0)=y/x$,
where $x$ and $y$ are determined by
\[ \left(\begin{array}{cc}-r_1&1\\-1&0\end{array}\right)
   \cdots
   \left(\begin{array}{cc}-r_n&1\\-1&0\end{array}\right)
   \left(\begin{array}{c}x\\y\end{array}\right)
   =
   \left(\begin{array}{c}-1\\1\end{array}\right) .\]
In analogy with the considerations above one finds
$y/x=[r_n,\ldots, r_2,r_1+1]$. Notice that
\[ [r_n,\ldots, r_2,r_1+1]=[r_n,\ldots ,r_{k+1},r_k+1],\]
where $k$ is the smallest index for which $r_k<-2$ (if all $r_i$ are
equal to~$-2$, then the continued fraction equals~$-1$).
As proved by Honda~\cite[Thm.~2.3]{hond00},
cf.~\cite{giro00}, the number of tight contact structures on the
solid torus with this convex boundary is equal to
\[ |(r_n+1)\cdots (r_1+1)| \]
(beware that Honda's notation differs slightly from ours).
We claim that this is exactly the number of tight contact structures
we can produce by choosing our surgery curves $K_1,\ldots ,K_n$
suitably. Indeed, we had fixed $tb(K_i)=r_i+1\leq -1$ (where we regard
$K_i$ as a Legendrian knot in $N_{i-1}'\equiv N_0\subset (S^3,\xi_0)$),
and then by
remark~(ii) in the preceding section we have a choice of $|r_i+1|$
different rotation numbers. One now argues as in the proof
of~\cite[Prop.~4.18]{hond00} (building on a result of Lisca
and Mati\'c) that the resulting contact structures are pairwise distinct.

This completes the proof of Proposition~\ref{prop:r<0}.

\vspace{2mm}

As a direct consequence of this proof we see that our problem of
replacing integer surgeries by $(\pm 1)$--surgeries has a simple
solution for negative surgery coefficients.

\begin{cor}
Any contact $n$--surgery with $n\in{\mathbf Z}^-$ can be replaced by
a suitable $(-1)$--surgery.
\end{cor}

Next we deal with positive integer surgery coefficients. Again we begin by
proving a more general statement about surgery with positive
rational coefficient.

\begin{prop}
\label{prop:r>0}
Any contact $r$--surgery with $r>0$ can be obtained by a contact
$(1/k)$--surgery with some positive integer~$k$, followed by a
suitable contact $r'$--surgery with~$r'<0$.
\end{prop}

\noindent {\bf Proof.}
Write $r=p/q$ with $p,q$ coprime positive integers. Choose integers
$p',q'$ such that $pq'-qp'=1$. Then, topologically, $r$--surgery
along a Legendrian knot $K$ is defined by cutting out a tubular
neighbourhood $\nu K$ and gluing back a solid torus $N_0$ with the
attaching map
\[ \mu_0\longmapsto p\mu +q\lambda ,\;\; \lambda_0\longmapsto p'\mu +
q'\lambda ,\]
where the notation is as in the proof of Proposition~\ref{prop:r<0};
in particular, $\lambda$ denotes the longitude determined by the
contact framing of~$K$.

This means that $-p'\mu_0+p\lambda_0$ is glued to the dividing curve~$\lambda$.
So the possible contact $r$--surgeries are determined by the tight
contact structures on $N_0$ with convex boundary satisfying $\#\Gamma_{\partial
N_0}=2$ and $s(\partial N_0)=-p/p'$.

Now choose a positive integer $k$ such that $q-kp<0$. We can perform
a $(1/k)$--surgery on $K$ by cutting out $\nu K$ and gluing in~$N_0$,
with gluing map
\[ \mu_0\longmapsto \mu +k\lambda ,\;\; \lambda_0\longmapsto\lambda .\]
The unique contact $(1/k)$--surgery is defined by the tight contact
structure on $N_0$ with convex boundary satisfying $\#\Gamma_{\partial
N_0}=2$ and $s(\partial N_0)=\infty$ (in terms of $(\mu_0,\lambda_0)$),
in other words, by taking $N_0$ to be the standard neighbourhood of
a Legendrian knot.

Next perform $r'$--surgery along the spine of~$N_0$,
with $r'=p/(q-kp)$. That is,
cut out a tubular neighbourhood $N_1$ of the spine of $N_0$ and glue back a
solid torus $N_1'$ with the attaching map
\[ \mu_1\longmapsto p\mu_0+(q-kp)\lambda_0,\;\;
\lambda_1\longmapsto p'\mu_0+(q'-kp')\lambda_0.\]
Since
\[ \left(\begin{array}{cc}1&0\\k&1\end{array}\right)
   \left(\begin{array}{cc}p&p'\\q-kp&q'-kp'\end{array}\right)
=
   \left(\begin{array}{cc}p&p'\\q&q'\end{array}\right) ,\]
the net result of these two surgeries is the desired $r$--surgery.

We see that $-p'\mu_1+p\lambda_1$ is glued to~$\lambda_0$, so again the
possible contact surgeries are the ones corresponding to tight contact
structures on the solid torus with boundary slope $-p/p'$. Moreover,
since rational contact surgery with negative coefficient can be effected
by a sequence of contact $(-1)$--surgeries, i.e.\ symplectic handlebody
surgeries, the result of performing a contact $r'$--surgery on $N_0$
(which can be realised on $N_0\subset (S^3,\xi_0)$, for instance) will
be a tight contact structure on all of~$N_0$. So the different choices
correspond exactly to the different contact $r$--surgeries
along~$K$.

\vspace{2mm}

Again we can specialise this to contact surgeries with positive integer
coefficients.

\begin{cor}
Any contact $n$--surgery with $n\in{\mathbf Z}^+$, $n\geq 2$, is
equivalent to a contact $(+1)$--surgery followed by a contact
$(n/(1-n))$--surgery, and hence (by Proposition~\ref{prop:r<0})
to a contact $(+1)$--surgery
followed by a sequence of contact $(-1)$--surgeries.
\end{cor}

The only case left open is that of a contact $0$--surgery. We claim that
the result of any contact $0$--surgery is the same as that of a
suitable contact $(+1)$--surgery.
We defer the proof of this fact to the next section, because it
can best be dealt with in the context of Lutz twists.
Taking this statement about $0$--surgery for granted,
we have proved that by performing contact $(\pm 1)$--surgeries
along a Legendrian link in~$M$, we can obtain $S^3$ with some contact
structure.
\section{Lutz twists and contact surgeries}
\label{section:lutz}
In order to prove Theorem~\ref{thm:main} it remains to show that,
given any contact structure $\xi_0'$ on~$S^3$, one can find a Legendrian
link in $S^3$ such that contact $(\pm 1)$--surgery along this link
produces $(S^3,\xi_0)$. For this it is sufficient to observe that
$(S^3,\xi_0')$ can be obtained from $(S^3,\xi_0)$ by suitable
so-called {\it Lutz twists}, and that any Lutz twist is equivalent to
certain integer surgeries. We shall now elaborate on these points.

First recall the concept of a Lutz twist. Let $K$ be any knot in
a contact $3$--manifold $(M,\xi )$ that is transverse to~$\xi$. Then
$K$ has a tubular neighbourhood $S^1\times D_{r_0}^2\subset M$, where
$D^2_{r_0}$ denotes the open $2$--disc of radius $r_0$ and
$K\equiv S^1\times\{ 0\}$, on which $\xi$ is given by $d\theta +r^2
d\varphi$ in the coordinates $(\theta ,r,\varphi )$ used before,
cf.~\cite{mart71}. By a diffeomorphism of the
solid torus $S^1\times D_{r_0}^2$ given by Dehn twists along the
meridian (which amounts to a different choice
of longitude) and rescaling the $r$--coordinate we may assume
that $r_0>1$.

A {\it simple Lutz twist} along $K$ is defined by replacing $\xi$ on
$S^1\times D^2_{r_0}$ by the contact structure $\beta_0=0$, where
$\beta_0=h_1(r)\, d\theta +h_2(r)\, d\varphi$ is defined as
in Section~\ref{section:surgery}, subject to the following
conditions (for some small $\varepsilon >0$):
\begin{itemize}
\item $h_1(r)\equiv -1$ and $h_2(r)\equiv -r^2$ for $r<\varepsilon$,
\item $h_1(r)\equiv 1$ and $h_2(r)=r^2$ for $r>1-\varepsilon$,
\item $(h_1(r),h_2(r))$ does {\it not} complete a full turn around
$(0,0)$ as $r$ goes from $0$ to~$r_0$.
\end{itemize}
Such a simple Lutz twist does not change the topology of the underlying
manifold, but it does, in general,
change the homotopy type of the contact structure
(as a $2$--plane field).

A {\it full Lutz twist} is defined by the conditions
\begin{itemize}
\item $h_1(r)\equiv 1$ and $h_2(r)=r^2$ for $r<\varepsilon$ and $r>
1-\varepsilon$,
\item $(h_1(r),h_2(r))$ completes exactly one full twist around $(0,0)$.
\end{itemize}
A full Lutz twist does not change the homotopy type of the contact
structure as a $2$--plane field, nor the topology of the underlying
manifold.

Notice that the disc $\{\theta_0\}\times D_{r_0'}^2$, where $r_0'>0$
is chosen such that $h_2(r_0')=0$, is an overtwisted disc (after a small
perturbation). In the case of a simple Lutz twist, $r_0'$ is
unique; in the case of a full Lutz twist there are two choices.

Homotopy classes of (cooriented) $2$--plane fields on $S^3$ are
classified by the homotopy group
$\pi_3(S^2)\cong {\mathbf Z}$. As shown by Lutz~\cite{lutz70},
cf.~\cite{lutz77} (again, the full details are folklore), each of these
homotopy classes contains a contact structure, obtained by suitable
Lutz twists from $(S^3,\xi_0)$. Moreover, Eliashberg~\cite{elia92}
has given a complete classification (up to isotopy) of contact structures
on~$S^3$: Each homotopy class contains exactly one overtwisted contact
structure, and the only tight contact structure is the standard one~$\xi_0$.
The overtwisted contact structure that is homotopic to $\xi_0$ as a $2$--plane
field can be obtained from $\xi_0$ by a full Lutz twist. In conclusion,
every contact structure on $S^3$ can be obtained from $\xi_0$ by Lutz
twists.

\begin{prop}
\label{prop:lutz}
A simple Lutz twist is equivalent to two contact $(+1)$--surgeries.
\end{prop}

\noindent{\bf Proof.} (cf.~\cite[Prop.~4.3]{etho})
Consider $N_0=S^1\times D^2_{r_0}$ with $r_0>1$ and contact structure
$\zeta_0$ given by $h_1(r)\, d\theta +h_2(r)\, d\varphi =0$, with $h_1,
h_2$ as described above for a simple Lutz twist.

Let $r_1$ be the unique radius for which $h_1(r_1)=-h_2(r_1)>0$. Then
along the $2$--torus $T_{r_1}=\{ r=r_1\}$ the contact structure
$\zeta_0$ is given by $d\theta -d\varphi =0$, so the characteristic
foliation on $T_{r_1}$ is linear of slope $d\theta /d\varphi =1$. As
$r$ goes from $0$ to~$r_1$, the slope of the characteristic foliation
on $T_r$ decreases monotonically from $0$ to $-\infty$, and then further
from $\infty$ to~$1$. This shows that the restriction of $\zeta_0$ to
$N_1=S^1\times D^2_{r_1}$ is tight, for by Dehn twists along the meridian
and rescaling the $r$--coordinate we can embed this into the solid torus
with its standard tight contact structure $d\theta +r^2\, d\varphi =0$.

Now perturb $T_{r_1}$ into a convex torus of slope~$1$ as described at the
end of Section~\ref{section:surgery}. Then the restriction of $\zeta_0$
to the perturbed $N_1$ is the tight contact structure uniquely
determined by this convex boundary and may be thought of as a
standard neighbourhood of its Legendrian spine~$K_1$. Notice that the
contact framing of $K_1$ is given by $\lambda_0+\mu_0$.

As a next step, we perform contact $(-1)$--surgery along~$K_1$. That is,
we cut out $N_1$ and glue back a solid torus $N_1'$ with the attaching map
\[ \mu_1\longmapsto\mu_0-(\lambda_0+\mu_0)=-\lambda_0,\;\;
\lambda_1\longmapsto \mu_0.\]
This means that $\lambda_1-\mu_1$ is mapped to the dividing curve
$\lambda_0+\mu_0$, so $N_1'$ carries the unique tight contact structure
with convex boundary of slope~$-1$.

We now claim that the net result of performing a simple Lutz twist and the
described contact $(-1)$--surgery is a contact $(+1)$--surgery. To
prove this, perturb the $2$--torus $T_1=\{ r=1\}$ into a convex torus
with dividing set of slope~$-1$ with respect to $(\mu_0,\lambda_0)$ (and
two components, as usual). Consider the thickened torus $T^2\times
[r_1,1]$ with boundaries the perturbed $T_{r_1}$ and~$T_1$. The restriction
of $\zeta_0$ to this thickened torus is tight (after a diffeomorphism of
the thickened torus it is seen to embed into the standard tight
structure on $S^1\times D^2$), the boundaries are convex, and $\zeta_0$
is {\it minimally twisting} in the sense of~\cite[Section~2.2.1]{hond00},
that is, the contact structure twists minimally in radial direction
to satisfy the boundary conditions. Notice that with respect to
$(\mu_1,\lambda_1)$ the slope of $T_1$ is~$+1$, since
$\lambda_1+\mu_1$ is glued to $\mu_0-\lambda_0$. Moreover, with respect
to $(\mu_1,\lambda_1)$ the slope twists from $-1$ on $T_{r_1}$ to
$+1$ on $T_1$ by passing via $\pm\infty$ rather than~$0$ ($\lambda_1$ is
glued to~$\mu_0$).

Let $r_1'$ be the radius determined by $h_1(r_1')=h_2
(r_1')<0$. We then see that $T^2\times [r_1,1]$ is contactomorphic to
$T^2\times [r_1',r_1]$ (either equipped with the restriction of~$\zeta_0$,
and with perturbed boundaries). Furthermore,
by Proposition~\ref{prop:1/k} we know that $N_1'$ is contactomorphic to
the perturbed $S^1\times D^2_{r_1'}$.
We conclude that
$N_1'\cup T^2\times [r_1,1]$ is contactomorphic to the perturbed~$N_1$,
i.e.\ the solid torus with the essentially unique tight contact structure
with convex boundary of slope~$+1$.

The perturbed $T_1$ has slope $-1$ with respect to
$(\mu_0,\lambda_0)$, so it may be regarded as the boundary of a standard
tubular neighbourhood of a Legendrian knot with contact framing
$\lambda_0-\mu_0$. So the result of performing a Lutz twist and
a $(-1)$--surgery as described is to cut out this tubular neighbourhood,
glue in a solid torus $N_2'$ with attaching map that sends $\mu_2$
to $-\lambda_0=-\mu_0-(\lambda_0-\mu_0)$, and extend with the unique
tight contact structure on $N_2'$ determined by this gluing. This
amounts to a contact $(+1)$--surgery.

We conclude that this last contact $(+1)$--surgery, followed by the contact $(+1)$--surgery reversing the $(-1)$--surgery along~$K_1$, amounts
to a simple Lutz twist.

\vspace{2mm}

As a more or less immediate corollary of this proof we have the following
results. They complete the proof of Theorem~\ref{thm:main}.

\begin{prop}
{\rm (i)} A full Lutz twist is equivalent to four contact $(+1)$--surgeries.

{\rm (ii)} A contact $0$--surgery is equivalent to a contact
$(+1)$--surgery.
\end{prop}

\noindent {\bf Proof.}
Argue as in the proof of the preceding proposition. For (i), perform
a contact $(-1)$--surgery determined by the torus $T_{r_1}$, where
$h_1(r_1)=-h_2(r_1)<0$, then a contact $(-1)$ surgery determined
by $T_{r_2}$, where $h_1(r_2)=h_2(r_2)<0$. This reduces a full Lutz twist
to a simple Lutz twist.

For (ii), perform a contact $(-1)$--surgery as in the case of a simple
Lutz twist. This reduces the contact $0$--surgery to a trivial
surgery. (The relation between the two Le\-gendrian knots corresponding
to the contact $0$--surgery and $(+1)$--surgery,
respectively, is analogous to that in the case of contact
surgeries with different negative integer framings, as discussed
in the proof of Proposition~\ref{prop:r<0}.)

\vspace{5mm}

\noindent {\bf Acknowledgements.}

\vspace{1mm}

This research was carried out while
F.~D.\ was staying at the Mathematical Institute of Leiden University
(then the home of H.~G.),
supported by a post-doctoral fellowship from NWO (Netherlands Organisation
for Scientific Research). The research of F.~D.\ is partially supported
by the NSFC (National Science Foundation of China).



\noindent Mathematisches Institut, Universit\"at zu K\"oln,
Weyertal 86--90, 50931 K\"oln, Germany; e-mail: geiges@math.uni-koeln.de

\noindent Department of Mathematics, Peking University, Beijing 100871,
P.~R.~China, e-mail: dingfan@math.pku.edu.cn 
\end{document}